\documentclass[12pt]{article}
\usepackage{amsmath,amstext, verbatim,remreset}
\usepackage[T2A]{fontenc}
\usepackage[utf8]{inputenc}
\usepackage[russian]{babel}
\usepackage{amsfonts,amssymb}
\usepackage[left=15mm, top=20mm, right=15mm, bottom=20mm]{geometry}
\usepackage{cite}
\newcommand{\sgn}{{\rm sgn}}
\newcommand{\RR}{\mathbb R}
\newcommand{\NN}{\mathbb N}
\newcommand{\ZZ}{\mathbb Z}
\newcommand{\K}{{\bf k}}
\newcommand{\M}{{\bf M}}
\usepackage{mathrsfs}
\newtheorem{theorem}{Теорема}
\newtheorem{lemma}{Лема}
\newtheorem{definition}{Означення}
\newtheorem{problem}{Задача}
\newtheorem{remark}{Зауваження}
\begin{document}
\begin{center}
	\begin{Large}
	{\bf Задачі найкращого відновлення на класах, що задаються обмеженнями на декілька старших похідних функцій}
	
	\bigskip
	В. Ф. Бабенко, О. В. Коваленко
	\end{Large}
\end{center}

\section{Постановка задачі і огляд відомих результатів}
Нехай задано підмножину $\mathfrak{M}$  множини $C$ неперервних $2\pi$ -- періодичних функцій $x(t)$ і числа $0\leq u_1 < u_2<\ldots<u_{2n}< 2\pi$, $u:=(u_1,\dots,u_{2n})$. 

Довільну функцію $\Phi\colon\RR^{2n}\to\RR$ ми будемо називати методом відновлення значення функції $x\in\mathfrak{M}$ в точці $\tau\in [0,2\pi)$, а величину $$e(\mathfrak{M}, u, \Phi, \tau):= \sup\limits_{x\in\mathfrak{M}}|\,x(\tau)-\Phi(x(u_1),\dots,x(u_{2n}))|$$ --- похибкою відновлення методом $\Phi$. 

Задача найкращого відновлення значення функції $x\in\mathfrak{M}$ в точці $\tau$ за її значеннями в  точках $u$ формулюється наступним чином.
\begin{problem}\label{ORinterpolation}
	Знайти найкращу похибку відновлення $$E(\mathfrak{M},u,\tau):= \inf\limits_\Phi e(\mathfrak{M}, u, \Phi, \tau)$$ а також найкращий метод відновлення $\tilde{\Phi}$, на якому досягається найкраща похибка відновлення. 
\end{problem}
Якщо нас цікавить вся функція $x(t)$, a не тільки її значення в фіксованій точці $\tau$, то виникає задача найкращого відновлення функції $x\in\mathfrak{M}$ за її значеннями в точках $u$. 

Довільну функцію $\Psi\colon \RR^{2n}\to C$ ми будемо називати методом відновлення функції $x\in\mathfrak{M}$, а величину $$e(\mathfrak{M}, u, \Psi,\|\cdot\|):= \sup\limits_{x\in\mathfrak{M}}\|\,x-\Psi(x(u_1),\dots,x(u_{2n}))\|$$ --- похибкою відновлення методом $\Psi$ ($\|\cdot\|$ --- деяка норма в просторі $C$).  
\begin{problem}\label{ORapproximation}
	Знайти найкращу похибку відновлення $$E(\mathfrak{M},u,\|\cdot\|):= \inf\limits_\Psi e(\mathfrak{M}, u, \Psi,\|\cdot\|)$$ а також найкращий метод $\tilde{\Psi}$, на якому досягається найкраща похибка відновлення. 
\end{problem}
Крім того, цікавим є питання про найкраще розташування інформаційних вузлів.
\begin{problem}\label{ORbest_information_set}
	Знайти вектор $u^*$, на якому досягається $$E(\mathfrak{M},\|\cdot\|):= \inf\limits_u E(\mathfrak{M}, u,\|\cdot\|).$$
\end{problem}

Відмітимо, що ми ставимо задачі відновлення за інформаційними множинами, що містять парну кількість вузлів. Це природно в силу специфіки періодичних функцій, але не є обов'язковим.


Відомі наступні результати, що стосуються вказаних задач найкращого відновлення.
\begin{enumerate}
	\item $\mathfrak{M} = W^r_\infty$. Випадок, коли $u$ --- рівномірне розбиття, розглянуто в роботах В.~М.~Тихомірова~\cite{Tihomirov} (задача~\ref{ORapproximation} --- в метриці простору $C$) і 
	О.~А.~Женсикбаєва~\cite{Zhensikbaev} (задача~\ref{ORapproximation} --- в метриці простору $L_p$, $1\leq p<\infty$). В~\cite{Velikin77} В.~Л.~Велікін, зокрема, показав, що розв'язком задачі~\ref{ORbest_information_set} в
	метриках просторів $C$ і $L_1$ є рівномірне розбиття. 
	\item $\mathfrak{M} = W^rH^\omega = \{x\in C^r\colon \omega(x^{(r)},t)\leq \omega(t)\}$, де $\omega(t)$ --- заданий модуль неперервності, $r=0,1$ (див. М.~П.~Корнейчук~\cite{Kornejchuk}).
	\item $\mathfrak{M} =W^L_{\infty,\sigma}=\{f\in L^{r}_\infty\colon |Lx(t)|\leq \sigma(t), t\in[0,2\pi]\}$, де $L$ --- диференціальний оператор порядку $r$, $\sigma(t)$ --- додатна неперервна функція (див.
	 Б.~Д.~Боянов~\cite{Bojanov}).
	 
\noindent У випадку, коли $u$ --- рівномірне розбиття, розв'язок задачі~\ref{ORapproximation} для наступних класів функцій міститься в монографії~\cite{MLD}.
	\item $\mathfrak{M} = W^r_2$ в метриці простору $L_2$.
	\item  $\mathfrak{M} = W^r_p$ в метриці простору $L_1$.
\end{enumerate}

Інші результати, що стосуються задач оптимального відновлення а також подальші посилання можна знайти в статтях \cite{Pinkus,Micchelli76a,Ligun,Bojanov75,Micchelli76b,Tikhomirov87} а також монографіях \cite{MLD,Osipenko,ZhensikbaevBook,Traub,Micchelli77}.

Для $r,d\in \NN$, цілих чисел $0<k_1<\ldots<k_d \leq r$, $\K:=(k_1,k_2,\dots,k_d)$, і додатних чисел $M_{k_1},\dots,M_{k_d}$, $\M:=(M_{k_1}, M_{k_2},\dots,M_{k_d})$,  позначимо через $W^r_\infty(\M,\K)$ клас $2\pi$ -- періодичних функцій $x\in L^r_\infty$ таких, що $\left\|x^{(k_i)}\right\|_\infty\leq M_{k_i}$ для всіх $k=1,\dots,d$.

Ми будемо розглядати наступні частинні випадки класів $W^r_\infty(\M,\K)$:
\begin{enumerate}
\item $d=2$, $k_1 = r - 1$,   $k_2 = r$;
\item $d=2$, $k_1 = r - 2$,   $k_2 = r$;
\item $d=3$, $k_1 = r - 2$,   $k_2 = r - 1$, $k_3 = r$.
\end{enumerate}

У всіх випадках, що розглядаються, ми будемо вважати, що $M_r = 1$. Замість $W^r_\infty(\M,\K)$ будемо писати  $W^r_{r-1}(M)$ у випадку, коли $\M = (M,1)$, $\K = (r-1,r)$;  $W^r_{r-2}(M)$  у випадку, коли $\M = (M,1)$, $\K = (r-2,r)$ і $W^r_{r-1,r-2}(M,N)$ у випадку, коли $\M = (M,N,1)$, $\K = (r-2,r-1,r)$.

Через $\nu(f)$ будемо позначати число істотних змін знаку періодичної функції $f$ на періоді; через $\nu(f,[a,b])$ --- число істотних змін знаку функції $f$ на відрізку $[a,b]$.

Ціллю даної статті є розв'язок задач~\ref{ORinterpolation} -- \ref{ORbest_information_set} для класів $W^r_{r-1}(M)$, $W^r_{r-2}(M)$ і $W^r_{r-1,r-2}(M,N)$.


\section {$W^r_\infty(\M,\K)$-ідеальні сплайни з заданими нулями}

\begin{definition}\label{ORspline_r,r-1}
	Функцію $x(t)\in W^r_{r-1}(M)$ будемо називати $W^r_{r-1}(M)$ -- ідеальним сплайном з вузлами в точках $t_1< t_2<\ldots<t_{2n+1} := 2\pi + t_1$, якщо для кожного $k=1,\dots, 2n$ існують числа $\alpha_k\in[0,t_{k+1} - t_k]$ 
	і $\varepsilon_k=\pm 1$ такі, що $x^{(r)}(t) = \varepsilon_k$ на інтервалі $(t_k, t_k+\alpha_k)$, $x^{(r)}(t) = 0$ на інтервалі $( t_k+\alpha_k, t_{k+1})$ і $x^{(r-1)}(t) = \varepsilon_k M$ на інтервалі $( t_k+\alpha_k, 
	t_{k+1})$.
\end{definition}
\begin{definition}\label{ORspline_r,r-2}
	Функцію $x(t)\in W^r_{r-2}(M)$ будемо називати $W^r_{r-2}(M)$ -- ідеальним сплайном с вузлами в точках $ t_1< t_2<\ldots<t_{2n+1} := 2\pi + t_1$, якщо для кожного $k=1,\dots, 2n$ існують числа $\alpha_k\in[0,t_{k+1} - 
	t_k]$ і $\varepsilon_k=\pm 1$  такі, що  $x^{(r)}(t) = \varepsilon_k$ на інтервалі $(t_k, t_k+\frac{\alpha_k}2)$, $x^{(r)}(t) = -\varepsilon_k$ на інтервалі $( t_k+\frac{\alpha_k}2, t_{k} + \alpha_k)$, $x^{(r)}(t) = 0$ 
	на інтервалі $( t_k+\alpha_k, t_{k+1})$ і $x^{(r-2)}(t) = \varepsilon_k M$ на інтервалі $( t_k+\alpha_k, t_{k+1})$.
\end{definition}
\begin{definition}\label{ORspline_r,r-1,r-2}
	Функцію $x(t)\in W^r_{r-1,r-2}(M,N)$ назвемо $W^r_{r-1,r-2}(M,N)$ -- ідеальним сплайном с вузлами в точках $ t_1< t_2<\ldots<t_{2n+1} := 2\pi + t_1$, якщо для кожного $k=1,\dots, 2n$ існують числа 
	$\alpha_k\in[0,t_{k+1} - t_k]$, $\beta_k\in[0,\alpha_k]$ і $\varepsilon_k=\pm 1$  такі, що  $x^{(r)}(t) = \varepsilon_k$ на інтервалі $(t_k, t_k+\frac{\beta_k}2)$, $x^{(r)}(t) = 0$ на інтервалі $(t_k+\frac{\beta_k}2, 
	t_k+\alpha_k-\frac{\beta_k}2)$ , $x^{(r)}(t) = -\varepsilon_k$ на інтервалі $(  t_k+\alpha_k-\frac{\beta_k}2, t_k+\alpha_k)$, $x^{(r)}(t) = 0$ на інтервалі $( t_k+\alpha_k, t_{k+1})$, $x^{(r-1)}(t) = \varepsilon_k N$ на 
	інтервалі $(t_k+\frac{\beta_k}2, t_k+\alpha_k-\frac{\beta_k}2)$  і  $x^{(r-2)}(t) = \varepsilon_k M$ на інтервалі $( t_k+\alpha_k, t_{k+1})$.
\end{definition}

\begin{theorem}\label{ORth::splineExistance}
	Нехай задано натуральні числа $n\in\NN$ і числа $0\leq u_1 < u_2<\ldots < u_{2n}<2\pi$. Нехай $X$ позначає один з класів $W^r_{r-1}(M)$, $W^r_{r-2}(M)$ або $W^r_{r-1,r-2}(M,N)$. Існує $X$ -- ідеальний сплайн з $2n$ 
	вузлами, що має нулі в точках  $u_1, u_2,\dots, u_{2n}$.
\end{theorem}

Нехай $X = W^r_{r-1}(M)$.

В просторі $\RR^{2n}_1$ розглянемо сферу $S^{2n-1}$ радіуса $2\pi$, тобто $$S^{2n-1}=\left\{ \xi=(\xi_1,\xi_2,\dots,\xi_{2n}):\sum_{i=1}^{2n}|\xi_i|=2\pi\right\}.$$ Кожній точці $\xi\in S^{2n-1}$ поставимо у відповідність розбиття відрізку $[0,2\pi]$ точками $t_0:=0,t_1:=|\xi_1|,t_2:=|\xi_1|+|\xi_2|,\dots,t_{2n-1}:=\sum_{i=1}^{2n-1}|\xi_i|,t_{2n}:=\sum_{i=1}^{2n}|\xi_i|=2\pi$. 
 
Для кожного $k=1,\dots,2n$ покладемо $\phi_1(\xi,t)=\sgn\xi_k\cdot\min\left(t-t_{k-1}, M\right)$ на відрізку $\left[t_{k-1}, \frac{t_{k-1}+t_k}{2}\right]$; $\phi_1(\xi,t)=\sgn\xi_k\cdot\min\left(t_k - t, M\right)$ на  відрізку $\left[\frac{t_{k-1}+t_k}{2},t_{k}\right]$. Таким чином ми визначили неперервну функцію $\phi_1(\xi,t)$ на всьому відрізку $[0,2\pi]$, причому для всіх $k=0,1,\dots,2n$ мають місце рівності $\phi_1\left(\xi,t_k\right)=0$. 
   
Для $k = 2,\dots,r-1$ покладемо $$\phi_k(\xi,t) = \int\limits_0^t\phi_{k-1}(\xi,s)ds - \frac 1{2\pi}\int\limits_0^{2\pi}\int\limits_0^t\phi_{k-1}(\xi,s)dsdt.$$ Покладемо $\phi_r(\xi,t) = \int\limits_{u_1}^t\phi_{r-1}(\xi,s)ds$. Визначимо відображення $\eta\colon S^{2n-1}\to \RR^{2n}$ наступним чином. $$\eta(\xi):=\left(\int\limits_0^{2\pi}\phi_1(\xi,t)dt, \phi_r(\xi,u_2), \phi_r(\xi,u_3),\dots, \phi_r(\xi,u_{2n})\right).$$

В силу побудови функції $\phi_r(\xi,t)$ відображення $\eta$ неперервне і непарне. Згідно з теоремою Борсука (див.~\cite{borsuk}) існує нуль $\xi^*\in S^{2n-1}$ функції $\eta$. Це означає, що $\phi_r(\xi^*,u_k)=0$, $k=1,\dots,2n$ і $\int\limits_0^{2\pi}\phi_1(\xi^*,t)dt=0$. Остання рівність дозволяє нам періодично продовжити функцію $\phi_r(\xi^*,t)$ на всю вісь зі збереженням гладкості. В силу побудови функцій $\phi_r(\xi,t)$ функція $\phi_r(\xi^*,t)$ --- $W^r_{r-1}(M)$ -- ідеальний сплайн c вузлами в точках розбиття. Крім того, $\phi_r(\xi^*,t)$ має нулі в заданих точках $u_1,u_2,\dots,u_{2n}$. Таким чином $\phi_r(\xi^*,t)$ --- шуканий $W^r_{r-1}(M)$ -- ідеальний сплайн.  Теорему у випадку, коли $X = W^r_{r-1}(M)$, доведено. 
   
Нехай $X = W^r_{r-2}(M)$.
   
Нехай $\alpha,\beta\in\RR$, $a\in [-M,M]$. Для $b\in\left[0,\frac{|\beta|}2\right]$ визначимо функцію $$\psi_b(\alpha,\beta;t)\colon [\alpha,\alpha + |\beta|\,]\to\RR$$ наступним чином.  $\psi_b(\alpha,\beta;t) = t-\alpha$ на відрізку $[\alpha, \alpha + b]$, $\psi_b(\alpha,\beta;t) = \alpha+2b -t$ на відрізку $[\alpha + b, \alpha+2b]$, $\psi_b(\alpha,\beta;t) = 0$ на відрізку $[\alpha + 2b, \alpha+|\beta|]$. 

Покладемо $$B:=\max\left\{b\in\left[0,\frac{|\beta|}2\right]\colon \left|a+\sgn\beta\int\limits_\alpha^t\psi_b(\alpha,\beta;s)ds\right|\leq M \forall t\in [\alpha,\alpha + |\beta|\,]\right\}$$
і $\psi(\alpha,\beta,a;t):=a+\sgn\beta\int\limits_\alpha^t\psi_B(\alpha,\beta;s)ds$.
Відмітимо, що справедливі наступні рівності:
\begin{equation}\label{ORthsplineExistance.1}
\psi'(\alpha,\beta,a;\alpha)=\psi'(\alpha,\beta,a;\alpha+|\beta|)=0.
\end{equation}

Для натурального числа $m\in\NN$ і  вектора $\zeta=(\zeta_1,\dots,\zeta_m)\in\RR^m$ визначимо функцію $$\psi(\alpha,\zeta,a;t)\colon \left[\alpha,\alpha + \sum\limits_{k=1}^m|\zeta_k|\,\right]\to\RR$$ наступним чином. На відрізку $[\alpha,\alpha + |\zeta_1|\,]$ $\psi(\alpha,\zeta,a;t)=\psi(\alpha, \zeta_1,a;t)$; на відрізку $[\alpha + |\zeta_1|,\alpha + |\zeta_1| + |\zeta_2|\,]$ $\psi(\alpha,\zeta,a;t)=\psi(\alpha + |\zeta_1|, \zeta_2,\psi(\alpha, \zeta,a;\alpha + |\zeta_1|);t)$; і так далі, на відрізку $\left[\alpha+\sum\limits_{k=1}^{m-1}|\zeta_k|,\alpha + \sum\limits_{k=1}^m|\zeta_k|\,\right]$ $$\psi(\alpha,\zeta,a;t)= \psi\left(\alpha+\sum\limits_{k=1}^{m-1}|\zeta_k|,\zeta_m,\psi\left(\alpha, \zeta,a; \alpha +\sum\limits_{k=1}^{m-1}|\zeta_k|\right);t\right).$$

Для вектора $\xi=(\xi_1,\dots,\xi_{2n})\in S^{2n-1}$ і $a\in [-M,M]$  покладемо $\psi(\xi,a;t) =\psi(0,\xi,a;t)$ (функція $\psi(\xi,a;t)$ визначена на відрізку $[0,2\pi]$).

Нехай $a\in (M,M + \pi]$ і  $\xi=(\xi_1,\dots,\xi_{2n})\in S^{2n-1}$.  Нехай, крім того, індекс $i\in\{0,1,\dots,2n-1\}$ такий, що $a-M\in\left[\sum\limits_{k=0}^i|\xi_k|,\sum\limits_{k=0}^{i+1}|\xi_k|\right)$ ($\xi_0:=0$). 

Визначимо функцію $\psi(\xi,a;t)\colon [0,2\pi]\to\RR$ наступним чином. Покладемо $\psi(\xi,a;t)=M$ на відрізку $[0,a-M]$ і  $\psi(\xi,a;t)=\psi(a-M,\zeta,M;t)$ на відрізку $[a-M,2\pi]$, де $$\zeta = \left(\sgn\xi_{i+1}\cdot\left(\sum\limits_{k=0}^{i+1}|\xi_k|-a+M\right),\xi_{i+2},\xi_{i+3},\dots,\xi_{2n}\right)\in\RR^{2n-i}.$$
Аналогічно визначимо функцію $\psi(\xi,a;t)\colon [0,2\pi]\to\RR$ для $\xi=(\xi_1,\dots,\xi_{2n})\in S^{2n-1}$ і $a\in [-M-\pi,-M)$: нехай індекс $i\in\{0,1,\dots,2n-1\}$ такий, що $-a-M\in\left[\sum\limits_{k=0}^i|\xi_k|,\sum\limits_{k=0}^{i+1}|\xi_k|\right)$ ($\xi_0:=0$). Покладемо $\psi(\xi,a;t)=-M$ на відрізку $[0,-a-M]$ і $\psi(\xi,a;t)=\psi(-a-M,\zeta,-M;t)$ на відрізку $[-a-M,2\pi]$, де $$\zeta = \left(\sgn\xi_{i+1}\cdot\left(\sum\limits_{k=0}^{i+1}|\xi_k|+a+M\right),\xi_{i+2},\xi_{i+3},\dots,\xi_{2n}\right)\in\RR^{2n-i}.$$
 
Таким чином для всіх $a\in [-M-\pi,M + \pi]$ і  $\xi\in S^{2n-1}$ ми визначили функцію $\psi(\xi,a;t)\colon [0,2\pi]\to\RR$. За побудовою функція $\psi(\xi,a;t)$ неперервна на $[0,2\pi]$. Крім того в силу~\eqref{ORthsplineExistance.1} вона неперервно диференційовна на $[0,2\pi]$. $|\psi(\xi,a;t)|\leq M$ при всіх $t\in [0,2\pi]$. Якщо $\xi\in S^{2n-1}$ і $-M-\pi\leq a_1<a_2\leq M+\pi$, то для всіх $t\in[0,2\pi]$ справедлива нерівність $\psi(\xi,a_1;t)\leq \psi(\xi,a_2;t)$. Відмітимо що, функції $\psi(\xi,a;t)$ неперервно залежать від параметрів $a$ і $\xi$. 
 
За побудовою для всіх $\xi\in S^{2n-1}$ $\int\limits_0^{2\pi}\psi(\xi,M+\pi;t)dt > 0$ і $\int\limits_0^{2\pi}\psi(\xi,-M-\pi;t)dt < 0$. Тому існує число $A=A(\xi)\in(-M-\pi,M+\pi)$ таке, що 
\begin{equation}\label{ORthsplineExistance.2}
	\int\limits_0^{2\pi}\psi(\xi,A;t)dt = 0.
\end{equation} 
Покладемо $\phi_2(\xi;t):=\psi(\xi,A;t)$. Для $k = 3,\dots,r-1$ покладемо $$\phi_k(\xi,t) = \int\limits_0^t\phi_{k-1}(\xi,s)ds - \frac 1{2\pi}\int\limits_0^{2\pi}\int\limits_0^t\phi_{k-1}(\xi,s)dsdt.$$ Покладемо $\phi_r(\xi,t) = \int\limits_{u_1}^t\phi_{r-1}(\xi,s)ds$. 
 
Відмітимо, що для $k=1,2,\dots,r-2$ і $k=r$ $\int\limits_0^{2\pi}\phi^{(k)}_r(\xi;t)dt = 0$.
 
Визначимо відображення $\eta\colon S^{2n-1}\to \RR^{2n}$ наступним чином. $$\eta(\xi):=\left(\int\limits_0^{2\pi}\phi^{(r-1)}_r(\xi,t)dt, \phi_r(\xi,u_2), \phi_r(\xi,u_3),\dots, \phi_r(\xi,u_{2n})\right).$$
В силу побудови функції $\phi_r(\xi,t)$ відображення $\eta$ неперервне і непарне. Згідно з теоремою Борсука існує нуль $\xi^*\in S^{2n-1}$ функції $\eta$. Це означає, що  $\phi_r(\xi^*,u_k)=0$, $k=1,\dots,2n$ і $\int\limits_0^{2\pi}\phi^{(r-1)}_r(\xi^*,t)dt=0$. Остання рівність дозволяє нам періодично продовжити функцію $\phi_r(\xi^*,t)$ на всю вісь зі збереженням гладкості. Крім того $\phi_r(\xi^*,t)$ дорівнює нулю в заданих точках $u_1,u_2,\dots,u_{2n}$.   
   
Якщо число $A=A(\xi^*)$, що фігурує у рівності~\eqref{ORthsplineExistance.2}, знаходиться на відрізку $[-M,M]$, то функція $\phi_r(\xi^*,t)$ --- $W^r_{r-2}(M)$ -- ідеальний сплайн з вузлами $0,t_1,\dots, t_{2n-1}$. Припустимо, що $A=A(\xi^*)\in  (M+\pi,M)$. Нехай $\xi^*=(\xi^*_1,\dots, \xi^*_{2n})$. Оскільки функція $\phi_r(\xi^*,t)$ має $2n$ нулів, то функція $\phi_r^{(r-2)}(\xi^*,t)$ має $2n$ змін знаку, а отже $A-M\in[0,|\xi^*_1|\,]$. Тоді функція $\phi_r(\xi^*,t)$ --- $W^r_{r-2}(M)$ -- ідеальний сплайн c вузлами в точках $t_1 - A + M,t_2,t_3,\dots, t_{2n-1}+A-M$. Аналогічно доводиться, що $\phi_r(\xi^*,t)$ --- $W^r_{r-2}(M)$ -- ідеальний сплайн і у випадку, коли $A=A(\xi^*)\in  (-M-\pi,-M)$.
   
Таким чином $\phi_r(\xi^*,t)$ --- шуканий $W^r_{r-2}(M)$ -- ідеальний сплайн. Теорему у випадку, коли $X = W^r_{r-2}(M)$, доведено.
   
Доведення теореми у випадку, коли $X = W^r_{r-1,r-2}(M,N)$ аналогічно доведенню у випадку $X = W^r_{r-2}(M)$. Треба лише замість функції $\psi_b(\alpha,\beta;t)$ розглядати функцію $\min\{\psi_b(\alpha,\beta;t),N\}$.
Теорему доведено.   

\begin{remark}
	Нехай задано числа $n\in\NN$ і  $0\leq u_1 < u_2<\ldots < u_{2n}<2\pi$. Нехай $X$  позначає один з класів $W^r_{r-1}(M)$, $W^r_{r-2}(M)$ або $W^r_{r-1,r-2}(M,N)$. Через $\phi(u;t) = \phi(X,u;t)$, 
	$u:=(u_1,\dots,u_{2n})$ будемо позначати $X$ -- ідеальний сплайн з $2n$ вузлами, що дорівнює нулю в точках $u_1,u_2,\dots,u_{2n}$. 
\end{remark}
 
\section{Інтерполяційні сплайни}
Нехай $X$ позначає один з класів $W^r_{r-1}(M)$, $W^r_{r-2}(M)$ або $W^r_{r-1,r-2}(M,N)$. 
\begin{definition} 
	Розбиття $$\Delta_{2n} \colon t_1<t_2<\ldots<t_{2n+1}=2\pi + t_1$$ будемо називати $X$ -- нормальним, якщо існує $X$ -- ідеальний сплайн $\phi$ з вузлами в точках розбиття $\Delta_{2n}$, що має $2n$ нулів на періоді.
\end{definition}
Нехай задано $W^r_{r-1}(M)$ -- нормальне розбиття $\Delta_{2n}$. Нехай $\phi(t)$ --- $W^r_{r-1}(M)$ -- ідеальний сплайн c вузлами в точках розбиття $\Delta_{2n}$, що має $2n$ нулів. Нехай, як і у  визначенні~\ref{ORspline_r,r-1}  $W^r_{r-1}(M)$ -- ідеальних сплайнів, числа $\alpha_k\in[0,t_{k+1} - t_k]$ такі, що $\phi^{(r)}(t) =\pm 1$ на інтервалі $(t_k, t_k+\alpha_k)$ і $\phi^{(r)}(t) = 0$ на інтервалі $( t_k+\alpha_k, t_{k+1})$, $k = 1,\dots, 2n$. 
 
\begin{definition}
  Функцію $s(t)\in L^{r}_\infty$ будемо називати $W^r_{r-1}(M)$ --- сплайном з вузлами в точках розбиття $\Delta_{2n}$, якщо $s^{(r)}(t) = c_k$ на інтервалі $(t_k, t_k+\alpha_k)$, $s^{(r)}(t) = 0$ на інтервалі $( t_k+\alpha_k, t_{k+1})$, $k = 1,\dots, 2n$, 
  $c_1,\dots, c_{2n}\in\RR$.
\end{definition}
Нехай задано $W^r_{r-2}(M)$ --- нормальне розбиття $\Delta_{2n}$. Нехай $\phi(t)$ --- $W^r_{r-2}(M)$ -- ідеальний сплайн з вузлами в точках розбиття $\Delta_{2n}$, що має $2n$ нулів. Нехай, як і у визначенні~\ref{ORspline_r,r-2}  $W^r_{r-2}(M)$ -- ідеальних сплайнів, числа $\alpha_k\in[0,t_{k+1} - t_k]$ такі, що   $\phi^{(r)}(t) = \pm 1$ на інтервалі $(t_k, t_k+\frac{\alpha_k}2)$, $\phi^{(r)}(t) = \mp 1$ на інтервалі $( t_k+\frac{\alpha_k}2, t_k+\alpha_k) $ і $\phi^{(r)}(t) = 0$ на інтервалі $( t_k+\alpha_k, t_{k+1})$, $k=1,\dots, 2n$. 
  
\begin{definition}
	Функцію $s(t)\in L^{r}_\infty$ будемо називати $W^r_{r-2}(M)$ --- сплайном з вузлами в точках розбиття $\Delta_{2n}$, якщо 
	\begin{enumerate}
		\item $s^{(r)}(t) = c_k^1$ на інтервалі $(t_{k}, t_{k}+\frac{\alpha_{k}}2)$ і $s^{(r)}(t) = c_k^2$ на інтервалі 
			 $( t_k+\frac{\alpha_k}2, t_k+\alpha_k)$, $c_k^1, c_k^2\in\RR$, $k=1,\dots, 2n$;
		\item $s^{(r)}(t) = 0$ на інтервалі $(t_k+\alpha_k, t_{k+1})$, $k=1,\dots, 2n$;
		\item $s^{(r-1)}(t_k) = 0$ для тих $k$, для яких $|\phi^{(r-2)}(t_k)|=M$;
		\item $c_k^2=c_{k+1}^1$ для тих $k$, для яких $|\phi^{(r-2)}(t_k)|<M$.
	\end{enumerate}		
\end{definition}
Нехай задано $W^{r}_{r-1,r-2}(M,N)$ --- нормальне розбиття $\Delta_{2n}$ і  $W^r_{r-1,r-2}(M,N)$ -- ідеальний сплайн $\phi(t)$ з вузлами в точках розбиття $\Delta_{2n}$, що має $2n$ нулів. Нехай, як і у визначенні~\ref{ORspline_r,r-1,r-2}  $W^r_{r-1,r-2}(M,N)$ -- ідеальних сплайнів, числа $\alpha_k\in[0,t_{k+1} - t_k]$, $\beta\in[0,\alpha_k]$ такі, що  $\phi^{(r)}(t) = \pm 1$ на множині $(t_k, t_k+\frac{\beta_k}2)$, $\phi^{(r)}(t) = \mp 1$, $(t_k+\alpha_k-\frac{\beta_k}2, t_k+\alpha_k)$ $\phi^{(r)}(t) = 0$ на множині $(t_k+\frac{\beta_k}2, t_k+\alpha_k-\frac{\beta_k}2)\bigcup ( t_k+\alpha_k, t_{k+1})$.
\begin{definition}
	Функцію $s(t)\in L^{r}_\infty$ будемо називати $W^r_{r-1,r-2}(M,N)$ --- сплайном з вузлами в точках розбиття $\Delta_{2n}$, якщо 
	\begin{enumerate}
		\item $s^{(r)}(t) = c_k^1$ на інтервалі $(t_{k}, t_{k}+\frac{\beta_{k}}2)$ і $s^{(r)}(t) = c_k^2$ на інтервалі 
			 $( t_k+\alpha_k - \frac{\beta_k}2, t_k+\alpha_k)$, $c_k^1, c_k^2\in\RR$, $k=1,\dots, 2n$;
		\item $s^{(r)}(t) = 0$ на інтервалах $( t_k+\alpha_k, t_{k+1})$ та  
			$( t_k+\frac{\beta_k}2, t_k-\frac{\beta_k}2+\alpha_k)$, $k=1,\dots, 2n$;
		\item $s^{(r-1)}(t_k) = 0$  для тих $k$, для яких $|\phi^{(r-2)}(t_k)|=M$;
		\item $c_k^2=c_{k+1}^1$ для тих $k$, для яких $|\phi^{(r-2)}(t_k)|<M$.
	\end{enumerate}		
\end{definition}

Позначимо через $S(X;\Delta_{2n})$ множину всіх $X$ -- сплайнів по розбиттю $\Delta_{2n}$. Відмітимо, що $S(X;\Delta_{2n})$ є $2n$ -- вимірним лінійним простором.
   
Справедливе наступне твердження.
\begin{lemma}\label{ORth::general_interpolation}
	Нехай $X$ позначає один з класів $W^r_{r-1,r-2}(M,N)$, $W^r_{r-1}(M)$ або $W^r_{r-2}(M)$ і задано $X$ -- нормальне розбиття $\Delta_{2n}$. 
	Нехай, крім того, точки $\tau_k$, $k=1,2,\dots, 2n$, задовольняють нерівностям $u_k<\tau_k<u_{k+1}$, $u_{2n+1} := u_1 + 2\pi$, де $u_k$ --- нулі 
	відповідного $X$ -- ідеального сплайна $\phi(X;t)$, $k=1,2,\dots, 2n$. Тоді для будь-якого набору чисел $y_1,\dots, y_{2n}$ існує єдиний 
	$X$ -- сплайн $s(t)\in S(X;\Delta_{2n})$, що задовольняє умовам $$s(\tau_k) = y_k,\,\,k=1,\dots,2n.$$
\end{lemma}
Для доведення леми~\ref{ORth::general_interpolation} нам знадобиться наступна лема.
\begin{lemma}\label{ORl::interpolation}
	Якщо виконані умови леми~\ref{ORth::general_interpolation} і для сплайна $s(t)\in S(X;\Delta_{2n})$  виконуються співвідношення
	\begin{equation}\label{ORl::interpolation.1}
		|s(\tau_k)| \leq |\phi(X;\tau_k)|,\,\, k=1,\dots, 2n,
	\end{equation}
	то $\|s^{(r)}\|_\infty\leq 1$.
\end{lemma}
Припустимо супротивне, нехай $\|s^{(r)}\|_\infty > 1$. Покладемо $s_*(t):=\varepsilon\frac{s(t)}{\|s^{(r)}\|_\infty}$, $\varepsilon=\pm 1$ (значення $\varepsilon$ ми виберемо пізніше).  Тоді
\begin{equation}\label{ORl::interpolation.2}
	\|s_*^{(r)}\|_\infty = 1.
\end{equation} 
 
В силу співвідношень~\eqref{ORl::interpolation.1} виконуються нерівності $|s_*(\tau_k)| < |\phi(X;\tau_k)|$, $k=1,\dots, 2n$. Це означає, що різниця $\delta(t) := s_*(t) - \phi(X;t)$ має не менше ніж $2n$ змін знаку, оскільки $\sgn \phi(X;\tau_k) = - \sgn\phi(X;\tau_{k+1})$, $k=1,\dots, 2n$. 

Для скорочення записів у доведенні цієї леми будемо писати $\phi(t)$ замість $\phi(X;t)$.
 
Нехай $X=W^r_{r-1}(M)$. Виберемо $\varepsilon$ так, щоб для деякого $k=1,\dots,2n$ майже всюди на інтервалі $t\in (t_k,t_{k+1})$ виконувалась рівність $s_*^{(r)}(t) = \phi^{(r)}(t)$. Це можливо в силу~\eqref{ORl::interpolation.2} і визначення $W^r_{r-1}(M)$ -- сплайнів. Функція $\delta^{(r)}$ не змінює знак всередині інтервалів $(t_j,t_{j+1})$ $j=1,\dots, 2n$ і рівна нулю на інтервалі $(t_k,t_{k+1})$. Це означає, що функція $\delta^{(r)}$ має не більш ніж $2n-2$ змін знаку, що неможливо.
 
Нехай  $X=W^r_{r-2}(M)$. Функція $\phi(t)$ має $2n$ нулів, тому для кожного $k=1,\dots, 2n$ похідна $\phi^{(r)}(t)$ має однаковий знак на інтервалах $( t_k+\frac{\alpha_k}2, t_k+\alpha_k)$ і $(t_{k+1}, t_{k+1}+\frac{\alpha_{k+1}}2)$. 

Якщо $|\phi^{(r-2)}(t_k)| < M$ для всіх $k=1,\dots, 2n$, то $\phi^{(r)}(t)$ майже всюди відмінна від нуля, а отже є кусково сталою з $2n$ змінами знаку. Тому в силу~\eqref{ORl::interpolation.2} ми можемо вибрати $\varepsilon=\pm 1$ так, що функція $\delta(t)$ має не більше $2n-2$ змін знаку, що не можливо.

Нехай натуральні числа $k_1<k_2$ такі, що $|\phi^{(r-2)}(t_{k_1})|= |\phi^{(r-2)}(t_{k_2})| = M$ і $|\phi^{(r-2)}(t_{k})|<M$, $k_1<k<k_2$. Тоді звуження функції $\phi^{(r-1)}(t)$ на відрізок $[t_{k_1}, t_{k_2-1} + \alpha_{k_2 -1}]$ є кусково лінійною функцією, на кожній ланці якої кутовий коефіцієнт дорівнює $1$ або $-1$, та яка змінює свій кутовий коефіцієнт $k_2-k_1$ разів. Відмітимо, що $\phi^{(r-1)}(t_{k_1}) = \phi^{(r-1)}(t_{k_2}) = 0$ і $s_*^{(r-1)}(t_{k_1}) = s_*^{(r-1)}(t_{k_2}) = 0$. Крім того, в силу~\eqref{ORl::interpolation.2} 
\begin{equation}\label{ORl::interpolation.21}
	\|s_*^{(r)}\|_{L_\infty(t_{k_1},t_{k_2})}\leq 1.
\end{equation}
Це означає, що якщо $S_1,\dots,S_m$ --- це такі множини додатної міри, що 
\begin{enumerate}
	\item $S_j\subset [t_{k_1}, t_{k_2}]$, $j=1,\dots, m$;
	\item Якщо $x\in S_j$, $y\in S_{j+1}$, то $x<y$, $j=1,\dots, m-1$;
	\item $\epsilon (-1)^j\delta^{(r-1)}(t) > 0$, $t\in S_j$, $\epsilon = \pm 1$, $j = 1,\dots, m$,
\end{enumerate}
то $m\leq k_2-k_1$. Більше того, якщо в~\eqref{ORl::interpolation.21} має місце рівність, то можна вибрати $\varepsilon =\pm 1$ так, щоб гарантувати нерівність $m\leq k_2-k_1 - 1$ --- для цього достатньо взяти $\varepsilon$ так, щоб на деякому проміжку кусково лінійні функції  $\phi^{(r-1)}$ та $s_*^{(r-1)}$ мали однаковий кутовий коефіцієнт. Це означає, що функція $\delta^{(r-1)}(t)$ (при відповідно вибраному значенні $\varepsilon$) має не більше ніж $2n-1$ змін знаку, що неможливо. Прийшли до суперечності

Випадок $X=W^r_{r-1,r-2}(M,N)$ доводиться аналогічно.
  
Лему доведено.
  
Повернемося до доведення леми~\ref{ORth::general_interpolation}. 

Доведемо, що тотожній нуль --- це єдиний сплайн з $S(X;\Delta_{2n})$, що задовольняє нульовим інтерполяційним умовам. Припустимо супротивне, нехай ненульовий сплайн $s(t)\in S(X;\Delta_{2n})$ такий, що $s(\tau_k) = 0$, $k=1,\dots, 2n$. Але тоді для будь-якого $\lambda\in\RR$ $\lambda s(\tau_k) = 0$, $k=1,\dots, 2n$, що суперечить лемі~\ref{ORl::interpolation}.

$S(X;\Delta_{2n})$ є лінійним простором розмірності $2n$. Нехай $\psi_1,\dots, \psi_{2n}$ --- деякий базис простору $S(X;\Delta_{2n})$. Тоді довільний $X$ -- сплайн $s(t)\in S(X;\Delta_{2n})$ можна представити у вигляді $s(t) =\sum\limits_{k=1}^{2n}c_k \psi_k(t)$.  Знаходження інтерполяційного сплайна зводиться до розв'язку системи $2n$ лінійних рівнянь $\sum\limits_{k=1}^{2n}c_k \psi_k(\tau_j) = y_j$, $j=1,\dots,2n$ відносно параметрів $c_k$, $k=1,\dots, 2n$, що визначають сплайн. Оскільки тотожній нуль --- це єдиний сплайн з $S(X;\Delta_{2n})$, що задовольняє нульовим інтерполяційним умовам, то однорідна система рівнянь має тільки нульовий розв'язок. Це означає, що при довільних числах $y_k$, $k=1,\dots, 2n$ відповідна неоднорідна система має єдиний розв'язок. Лему доведено.
 
\begin{definition}
	Для класу $X$ $2\pi$ -- періодичних диференційовних функцій покладемо $$X':=\{x'(t)\colon x(t)\in X\}.$$
\end{definition}
\begin{theorem}\label{ORth::function_interpolation}
	Нехай $X$ позначає один з класів $W^r_{r-1,r-2}(M,N)$, $W^r_{r-1}(M)$ або $W^r_{r-2}(M)$. Нехай задано функція $x(t)\in X$ і числа 
	$u_1<u_2<\ldots<u_{2n}<u_1+2\pi$, $u=(u_1,\dots, u_{2n})$. Нехай $\Delta_{2n} \colon t_1<t_2<\ldots<t_{2n}<2\pi + t_1$ --- вузли 
	$X$ -- ідеального сплайна $\phi(X;u)$, що дорівнює нулю в точках $u_1,\dots, u_{2n}$. Тоді існує єдиний $X'$ -- сплайн 
	$s(t)\in S(X',\Delta_{2n})$ такий, що $s(u_k) = x(u_k)$, $k=1,\dots,2n$.
\end{theorem}
 
Функція $\phi'(X;u)$ є $X'$ -- ідеальним сплайном, що має $2n$ нулів. Тому розбиття $\Delta_{2n}$ є $X'$ -- нормальним. Крім того, інтерполяційні вузли $u_1,\dots, u_{2n}$ є нулями функції $\phi(X;u)$, а отже в силу теореми Ролля розташовані між нулями $X'$ -- ідеального сплайна $\phi'(X;u)$. Це означає, що виконуються умови леми~\ref{ORth::general_interpolation}, з якої ми одразу отримуємо існування (і єдиність) шуканого сплайна $s(t)\in S(X',\Delta_{2n})$, що інтерполює функцію $x(t)$ в точках $u_1,\dots,u_{2n}$. Теорему доведено.
%
%
%
\begin{remark}
	В умовах теореми~\ref{ORth::function_interpolation} сплайн $s(t)\in S(X',\Delta_{2n})$, який інтерполює функцію $x$ в точках $u$,  ми будемо 	позначати через $s(x,u;t)=s(X',x,u;t)$. Крім того, сплайн $s(t)\in S(X',\Delta_{2n})$, що інтерполює значення $v=(v_1,\dots, v_{2n})$ в точках 
	$u$, ми будемо позначати $s(X',u,v;t)$.
\end{remark}
  
\section{Екстремальні властивості ідеальних сплайнів}
$X$ -- ідеальні сплайни мають наступну екстремальну властивість. 
 
\begin{theorem}\label{ORth::extremalProperty}
	Нехай $X$ позначає один з класів $W^r_{r-1,r-2}(M,N)$, $W^r_{r-1}(M)$ або $W^r_{r-2}(M)$. Нехай задано натуральне число $n\in\NN$ і числа 
	$0\leq u_1 < u_2<\ldots < u_{2n}<2\pi$, $u=(u_1,\dots, u_{2n})$. Тоді для всіх $x\in X$ і $t\in[0,2\pi)$ $$|x(t) - s(X',x,u;t)|\leq |\phi(X,u;t)|.$$
\end{theorem}
 
Припустимо супротивне, нехай існує точка $t^*\in[0,2\pi)$ така, що $|x(t^*)- s(X',x,u;t^*)| > |\phi(X,u;t^*)|$. Тоді існує число $\lambda\in(-1,1)$ таке, що $\lambda (x(t^*) - s(X',x,u;t^*)) = \phi(X,u;t^*)$. 

Покладемо $\delta(t):=\lambda  x(t) - \lambda s(X',x,u;t) - \phi(X,u;t)$. Функція $\delta(t)$ має $2n+1$ нуль: точки $u_1,\dots,u_{2n}$ і точка $t^*$.

Для скорочення позначень у доведенні цієї теореми замість $\phi(X,u;t)$ будемо писати просто $\phi(t)$, замість $s(X',x,u;t)$ --- просто $s(t)$. 

Нехай $X = W^r_{r-1}(M)$.

В силу теореми Ролля функція $\delta^{(r-1)}(t)$ має не менше ніж $2n+1$ змін знаку. 

Нехай $k\in \{1,2,\dots, 2n\}$. Ми можемо вважати, що $\phi^{(r)}(t) = 1$ на інтервалі $(t_k,t_k+\alpha_k)$, $\alpha_k\in[0,t_{k+1} - t_k]$ і $\phi^{(r)}(t) = 0$ на інтервалі $(t_k+\alpha_k, t_{k+1})$. Тоді функція $\phi^{(r-1)}(t) +\lambda s^{(r-1)}(t)$ є лінійною з кутовим коефіцієнтом $1$ на інтервалі $(t_k,t_k+\alpha_k)$ і $\phi^{(r-1)}(t) + \lambda s^{(r-1)}(t) = M$ на інтервалі $(t_k+\alpha_k, t_{k+1})$. Враховуючи те, що $x(t)\in W^r_{r-1}(M)$, ми отримуємо, що $|\lambda x^{(r)}(t)| < 1$  і $|\lambda x^{(r-1)}(t)| < M$ для майже всіх $t$. Це означає, що $\delta^{(r-1)}(t) < 0$ на інтервалі $(t_k+\alpha_k, t_{k+1})$ і функція $\delta^{(r-1)}(t)$ має не більше ніж одну зміну знаку на інтервалі $(t_k,t_k+\alpha_k)$, причому якщо функція $\delta^{(r-1)}(t)$ змінює знак на інтервалі $(t_k,t_k+\alpha_k)$, то ця зміна знаку з ''плюс'' на ''мінус''. Таким чином на кожному проміжку $(t_k,t_{k+1})$ між вузлами сплайна $\phi(t)$ функція $\delta^{(r-1)}$ має не більше ніж одну зміну знака, і, крім того, якщо на кожному з проміжків $(t_k,t_{k+1})$ і $(t_{k+1},t_{k+2})$ функція $\delta^{(r-1)}$ має по одній зміні знаку, то ця функція не має зміни знаку в точці $t_{k+1}$. Отже $\nu(\delta^{(r-1)})\leq 2n$. Прийшли до суперечності. Теорему у випадку, коли $X = W^r_{r-1}(M)$, доведено.

Нехай $X = W^r_{r-2}(M)$.

Припустимо, що $\left|\phi^{(r-2)}(t)\right|< M$ на відрізку $[0,2\pi]$. 
Тоді для кожного $k=1,\dots, 2n$ функція $\phi^{(r-1)}(t) + \lambda s^{(r-1)}(t)$ є лінійною з кутовим коефіцієнтом $\pm 1$ на інтервалі $\left(\frac {t_{k+1} + t_{k}} 2, \frac {t_{k+2} + t_{k+1}} 2\right)$, причому на сусідніх інтервалах вказаного виду кутові коефіцієнти мають протилежні знаки. Це означає, що на кожному інтервалі $\left(\frac {t_{k+1} + t_{k}} 2, \frac {t_{k+2} + t_{k+1}} 2\right)$ функція $\delta^{(r-1)}(t)$ має не більш ніж одну зміну знаку, причому з ''плюс'' на ''мінус'', якщо кутовий коефіцієнт дорівнює $1$ і з  ''мінус'' на ''плюс'', якщо кутовий коефіцієнт дорівнює $-1$. Звідси ми отримуємо, що $\nu(\delta^{(r-1)})\leq 2n$, що неможливо. 

Таким чином $\max\limits_{t\in [0,2\pi]}|\phi^{(r-2)}(t)| = M$. З теореми Ролля слідує, що функція $\delta^{(r-2)}(t)$ має не менше $2n+1$ змін знаку.

Нехай відрізок $[\alpha,\beta]\subset\RR$ такий, що 
\begin{equation}\label{ORthextremalProperty.1.1}
	\left|\phi^{(r-2)}(\alpha)\right|=\left|\phi^{(r-2)}(\beta)\right|=M
\end{equation}
 і 
\begin{equation}\label{ORthextremalProperty.1.2}
	\left|\phi^{(r-2)}(t)\right|<M,\,\,t\in (\alpha,\beta).
\end{equation}
Можемо вважати, що $\alpha = t_0 < t_1<\ldots<t_{l-1} < \beta$  --- всі вузли $W^r_{r-2}(M)$ -- ідеального сплайна $\phi(t)$ на інтервалі $[\alpha,\beta)$ (у випадку необхідності можемо перенумерувати вузли). Покладемо $t_l :=\beta$ (число $t_l$ може бути вузлом $\phi(t)$, а може їм не бути). Покажемо, що $\nu\left(\delta^{(r-2)},[\alpha,\beta]\right)\leq l$. Відмітимо, що в силу~\eqref{ORthextremalProperty.1.2} і означення $W^r_{r-2}(M)$ -- ідеальних сплайнів функція $\phi^{(r)}$ відмінна від нуля майже всюди на $(\alpha,\beta)$.
 
Без зменшення загальності можемо вважати, що 
\begin{equation}\label{ORthextremalProperty.2}
	\phi^{(r-2)}(\alpha) = -M.
\end{equation}
Покладемо $q_0 :=t_0$, $q_k:=\frac{t_{k-1}+t_k}2$, $k=1,\dots,l$, $q_{l+1} = t_{l}$, $\psi(t) := \phi(t) + \lambda s(t)$.

В силу побудови функцій $\phi$ і $s$ і рівності~\eqref{ORthextremalProperty.2} на інтервалі $(q_k,q_{k+1})$ 
\begin{equation}\label{ORthextremalProperty.3}
	\psi^{(r)}(t)=(-1)^k,\,k=0,1,\dots,l.
\end{equation} 
Оскільки 
\begin{equation}\label{ORthextremalProperty.4}
	\|\lambda x^{(r)}\|_\infty< 1,
\end{equation}
то на кожному з проміжків $(q_k,q_{k+1})$ функція $\delta^{(r-1)}$ може змінювати знак не більше ніж один раз, причому з ''плюс'' на ''мінус'' при парних $k$ і з ''мінус'' на ''плюс'' при непарних $k$, $k=0,\dots,l$. Таким чином $\nu\left(\delta^{(r-1)},[\alpha,\beta]\right)\leq l+1
$. Це означає, що $\nu\left(\delta^{(r-2)},[\alpha,\beta]\right)\leq l+2$. 

Припустимо, що $\nu\left(\delta^{(r-2)},[\alpha,\beta]\right)=l+2$. Тоді існують точки $\alpha\leq r_1<r_2<\ldots<r_{l+2}\leq\beta$ такі, що  $\delta^{(r-2)}(r_1)=\delta^{(r-2)}(r_2)=\ldots=\delta^{(r-2)}(r_{l+2})=0$  і (в силу~\eqref{ORthextremalProperty.2})  $(-1)^k\delta^{(r-1)}(r_k)>0$, $k=1 ,\dots,l+2$. Це означає, що $\nu\left(\delta^{(r-1)},[\alpha,\beta]\right)= l + 1$ і перша зміна знаку $\delta^{(r-1)}$ на відрізку $[\alpha,\beta]$ відбувається з ''мінус'' на ''плюс'', що неможливо. Таким чином $\nu\left(\delta^{(r-2)},[\alpha,\beta]\right)\leq l+1$.

Припустимо, що $\nu\left(\delta^{(r-2)},[\alpha,\beta]\right)= l+1$. Враховуючи~\eqref{ORthextremalProperty.2} отримаємо, що 
\begin{equation}\label{ORthextremalProperty.5}
	(-1)^{l}\delta^{(r-2)}(\beta)\leq 0.
\end{equation}

Крім того, функція $\delta^{(r-1)}(t)$ змінює знак не менше $l$ раз, причому перша зміна знаку з ''мінус'' на ''плюс''. Це означає, що існує точка $\tau\in(q_l,q_{l+1})$ така, що $(-1)^l\delta^{(r-1)}(\tau)<0$. Оскільки виконуються співвідношення~\eqref{ORthextremalProperty.3} і~\eqref{ORthextremalProperty.4}, то  $(-1)^l\delta^{(r-1)}(\beta)<0$. Але звідси в силу~\eqref{ORthextremalProperty.1.1},~\eqref{ORthextremalProperty.2} і~\eqref{ORthextremalProperty.5} слідує, що $\|x^{(r-2)}\|_\infty>M$, що неможливо. Це означає, що $\nu\left(\delta^{(r-2)},[\alpha,\beta]\right)\leq l$, а отже $\nu\left(\delta^{(r-2)}\right)\leq 2n$. Прийшли до суперечності. Теорему у випадку, коли $X = W^r_{r-2}(M)$, доведено. 
 
Випадок, коли $X = W^r_{r-1,r-2}(M,N)$ можна довести аналогічно до випадку $X = W^r_{r-2}(M)$.
 
\begin{theorem}\label{ORth::interiorComperison}
	Нехай $X$ позначає один з класів  $W^r_{r-1,r-2}(M,N)$, $W^r_{r-1}(M)$ або $W^r_{r-2}(M)$. Нехай задано натуральне число $n\in\NN$ і числа 
	$0\leq u_1 < u_2<\ldots < u_{2n}<2\pi$, $u:=(u_1,\dots,u_{2n}), u^*:=\left(0,\frac{\pi}{n},\frac{2\pi}{n},\dots,\frac{(2n-1)\pi}{n}\right)$. 
	Тоді якщо $u\neq u^*$, то  $$\|\phi(X,u)\|_\infty > \|\phi(X,u^*)\|_\infty.$$
\end{theorem}
Доведемо теорему у випадку, коли $X=W^r_{r-2}(M)$. Інші випадки доводяться аналогічно. Для скорочення записів будемо замість $\phi(X,u)$ писати  $\phi(u)$.
 
Відмітимо, що функція $\phi(u^*)$ є $\frac{2\pi}{n}$ періодичною і, крім того, справедлива рівність $\phi(u^*;t) = -\phi\left(u^*;t+\frac\pi n\right)$. 
 
Припустимо супротивне. Нехай $u\neq u^*$ і 
\begin{equation}\label{ORthinteriorComperison.1}
	\|\phi(u)\|_\infty \leq \|\phi(u^*)\|_\infty.
\end{equation} 
Тоді для будь-якого $\gamma\in\RR$
\begin{equation}\label{ORthinteriorComperison.2}  
	\nu(\pm\phi(u;\gamma+\cdot)-\phi(u^*;\cdot))\geq 2n.
\end{equation}
Можливі два випадки.
   
{\it Випадок 1}. Нехай $\|\phi^{(r-2)}(u^*;\cdot)\|_\infty=M$. Без зменшення загальності можемо вважати, що функція $\phi^{(r-2)}(u^*;t)$ зростає на відрізку $[0,\frac\pi n-\alpha]$, постійна на $[\frac \pi n-\alpha,\frac \pi n]$, спадає на $[\frac \pi n,\frac {2\pi} n-\alpha]$ і постійна на $[\frac {2\pi} n-\alpha, \frac {2\pi} n]$, де $\alpha \geq 0$. За означенням $W^r_{r-2}(M)$ -- ідеального сплайна 
\begin{equation}\label{ORthinteriorComperison.1.1}
	|\phi^{(r-2)}(u^*;t)|=M,\, t\in \left[\frac \pi n-\alpha,\frac \pi n\right]\bigcup \left[\frac {2\pi} n-\alpha, \frac {2\pi} n\right]. 
\end{equation}  
Нехай $ t_1< t_2<\ldots<t_{2n+1} := 2\pi + t_1$ --- вузли  $W^r_{r-2}(M)$ -- ідеального сплайна $\phi(u;t)$. Нехай індекс  $i\in\{1,2,\dots,2n\}$ такий, що  $t_{i+1}-t_i>\frac{\pi} n$. Без зменшення загальності можемо вважати, що   $\phi^{(r-2)}(u;t)$ не спадає на проміжку $[t_i,t_{i+1}]$. Тоді існує число $\beta>\alpha$ таке, що $\phi^{(r-2)}(u;t)$ зростає на проміжку $[t_i,t_{i+1}-\beta]$ і $\phi^{(r-2)}(u;t)=M$ на проміжку $[t_{i+1}-\beta,t_{i+1}]$. Розглянемо функцію $\delta(t):=\phi^{(r-2)}(u;t - t_i) - \phi^{(r-2)}(u^*;t)$. З \eqref{ORthinteriorComperison.2} випливає, що $\nu(\delta)\geq 2n$. Відмітимо, що звуження функції $\phi(X,u^*)$ на відрізок $\left[\frac{k\pi}{n}, \frac{(k+1)\pi}{n}-\alpha\right]$, $k\in\ZZ$, співпадає зі звуженням на цей відрізок відповідно підібраного ідеального сплайна Ейлера. Тому з теореми порівняння Колмогорова (див.~\cite{Kolmogorov85}) випливає, що на відрізку $\left[\frac{k\pi}{n}, \frac{(k+1)\pi}{n}-\alpha\right]$ функція $\delta(t)$ може мати не більше ніж одну зміну знака (причому з  ''плюс'' на ''мінус'', якщо  $\phi^{(r-1)}(u^*;t)$ невід'ємна на цьому проміжку і з ''мінус'' на ''плюс'', якщо недодатна). Тепер, враховуючи~\eqref{ORthinteriorComperison.1.1} і симетрії функції $\phi(u^*;t)$, ми отримаємо, що на кожному з проміжків знакопостійності функції $\phi^{(r-1)}(u^*;t)$ функція $\delta(t)$ може мати не більше ніж одну зміну знака, причому з ''плюс'' на ''мінус'' на проміжках невід'ємності $\phi^{(r-1)}(u^*;t)$ і з ''мінус'' на ''плюс'' на проміжках недодатності $\phi^{(r-1)}(u^*;t)$. Але $[t_{i+1}-t_i-\beta,t_{i+1}-t_i]\supset [\frac \pi n-\alpha,\frac \pi n]$ і $\delta(t) \geq 0$ для всіх $t\in [t_{i+1}-t_i-\beta,t_{i+1}-t_i]$. Це означає, що функція $\delta(t)$ не має змін знаку на проміжку $[0,\frac\pi n]$. Але тоді $\nu(\delta)\leq 2n-1$. Прийшли до суперечності.
  
{\it Випадок 2}. Нехай $\|\phi^{(r-2)}(u^*;\cdot)\|_\infty<M$. Тоді $\phi(u^*;t)$ --- ідеальний сплайн Ейлера. З~\eqref{ORthinteriorComperison.1} і теореми порівняння Колмогорова випливає, що $$\|\phi^{(r-2)}(u;\cdot)\|_\infty\leq \|\phi^{(r-2)}(u^*;\cdot)\|_\infty<M.$$  Це означає, що $|\phi^{(r)}(u;t)|_\infty=1$ майже всюди на $[0,2\pi]$.

Нехай $ t_1< t_2<\ldots<t_{2n+1} := 2\pi + t_1$ --- вузли $W^r_{r-2}(M)$ -- ідеального сплайна $\phi(u;t)$. Тоді існує індекс $i\in\{1,2,\dots,2n\}$ такий, що $t_{i+1}-t_i>\frac{\pi} n$. Тоді можна вибрати зсув $\phi(u;\gamma+\cdot)$ функції $\phi(u;t)$ так, що $ \nu(\phi^{(r)}(u;\gamma+\cdot)-\phi^{(r)}(u^*;\cdot))\leq 2n-1$. Але це суперечить~\eqref{ORthinteriorComperison.2}. Теорему доведено.
\section{Задачі найкращого відновлення}
Для розв'язку задачі~\ref{ORinterpolation} нам знадобиться наступна лема.
\begin{lemma}\label{ORduality}
	Нехай $X$ позначає один з класів $W^r_{r-1,r-2}(M,N)$, $W^r_{r-1}(M)$ або $W^r_{r-2}(M)$. Нехай задано числа 
	$0\leq u_1 < u_2<\ldots<u_{2n}< 2\pi$, $u:=(u_1,\dots,u_{2n})$ і $\tau\in [0,2\pi)$.  Тоді  
	$$E(X,u,\tau)=\sup \left\{ x(\tau)\colon x\in X, x(u_k)=0, k=1,\dots, 2n\right\}.$$
\end{lemma} 
Відмітимо, що класи $W^r_{r-1}(M)$, $W^r_{r-2}(M)$ або $W^r_{r-1,r-2}(M,N)$ є опуклими і центрально-симетричними. Тому серед оптимальних методів відновлення існує лінійний (див.~\cite{smolyk}). Це означає, що $$E(X,u,\tau)=\inf\limits_{c_k}\sup\limits_{x\in X} \left(x(\tau)-\sum\limits_{k=1}^{2n}c_kx(u_k)\right),$$
де точна нижня межа береться по всім векторам $(c_1,\dots, c_{2n})\in\RR^{2n}$.

Клас $X$ можна представити у наступному вигляді: $X = \left\{x\in L^r_\infty\colon p(x)\leq 1\right\}$,  де $$p(x)= \max\left\{ \|x^{(r)}\|_\infty, \frac{\|x^{(r-1)}\|_\infty}M\right\}$$ у випадку, коли $X = W^r_{r-1}(M)$; $p(x)= \max\left\{ \|x^{(r)}\|_\infty, \frac{\|x^{(r-2)}\|_\infty}M\right\}$ у випадку, коли $X = W^r_{r-2}(M)$ і $$p(x)= \max\left\{ \|x^{(r)}\|_\infty, \frac{\|x^{(r-1)}\|_\infty}N, \frac{\|x^{(r-2)}\|_\infty}M\right\}$$  у випадку, коли $X = W^r_{r-1,r-2}(M,N)$. Відмітимо, що у всіх трьох випадках функція $p(x)$ додатно однорідна і напівадитивна.

Тепер, застосовуючи міркування аналогічні до доведення теореми~1.3.4 монографії~\cite{ExactConstants} (тільки потрібно застосувати більш загальне формулювання теореми Хана-Банаха, див., напр., глава~3 \S~2 в монографії~\cite{hahn-banach}), отримуємо справедливість твердження леми. Лему доведено.
  
З леми~\ref{ORduality} і теорем~\ref{ORth::function_interpolation} і~\ref{ORth::extremalProperty} отримуємо розв'язок задачі~\ref{ORinterpolation}.
\begin{theorem}\label{ORth::pr1_solution}
	Нехай $X$ позначає один з класів $W^r_{r-1,r-2}(M,N)$, $W^r_{r-1}(M)$ або $W^r_{r-2}(M)$. Нехай задано числа 
	$0\leq u_1 < u_2<\ldots<u_{2n}< 2\pi$, $u:=(u_1,\dots,u_{2n})$ і $\tau\in [0,2\pi)$.  Тоді $$E(X,u,\tau)=|\phi(X,u;\tau)|.$$
	Найкращим  методом відновлення у задачі~\ref{ORinterpolation} є метод $$\tilde{\Phi}(v_1,\dots,v_{2n}) = s(X',u,v;\tau).$$
\end{theorem}
  
З теорем~\ref{ORth::function_interpolation},~\ref{ORth::extremalProperty} і~\ref{ORth::pr1_solution} отримуємо розв'язок задачі~\ref{ORapproximation}.
   
\begin{theorem}\label{ORth::pr2_solution}
	Нехай $X$ позначає один з класів $W^r_{r-1,r-2}(M,N)$, $W^r_{r-1}(M)$ або $W^r_{r-2}(M)$. Нехай задано числа 
	$0\leq u_1 < u_2<\ldots<u_{2n}< 2\pi$, $u:=(u_1,\dots,u_{2n})$ і $p\in [1,\infty]$.  Тоді $$E(X,u,\|\cdot\|_p)=\|\phi(X,u;\cdot)\|_p.$$
	Найкращим методом відновлення у задачі~\ref{ORapproximation} є метод $$\tilde{\Psi}(v_1,\dots,v_{2n}) = s(X',u,v;t).$$
\end{theorem}
З теорем~\ref{ORth::interiorComperison} і~\ref{ORth::pr2_solution} отримуємо наступну теорему.
\begin{theorem}\label{ORth::pr3_solution}
	Нехай $X$ позначає один з класів $W^r_{r-1,r-2}(M,N)$, $W^r_{r-1}(M)$ або $W^r_{r-2}(M)$. Нехай 
	$u^*:=\left(0,\frac{\pi}{n},\frac{2\pi}{n},\dots,\frac{(2n-1)\pi}{n}\right)$. Тоді 
	$$E(X,\|\cdot\|_\infty):= e(X, u^*,\|\cdot\|_\infty),$$
	тобто найкращою інформаційною множиною (у випадку рівномірної метрики) є рівномірне розбиття відрізку $[0,2\pi]$.
\end{theorem}
\bibliographystyle{my_ugost2003s}
\bibliography{bibliography}

\begin{thebibliography}{10}
\def\selectlanguageifdefined#1{
\expandafter\ifx\csname date#1\endcsname\relax
\else\selectlanguage{#1}\fi}
\providecommand*{\href}[2]{{\small #2}}
\providecommand*{\url}[1]{{\small #1}}
\providecommand*{\BibUrl}[1]{\url{#1}}
\providecommand{\BibAnnote}[1]{}
\providecommand*{\BibEmph}[1]{#1}

\bibitem{Velikin77}
\selectlanguageifdefined{russian}
\BibEmph{Великин~В.~Л.} Оптимальная
  интерполяция периодических
  дифференцируемых функций с ограниченной
  старшей производной~/ В.~Л.~Великин~//
  \BibEmph{Матем. заметки}. ---
\newblock 1977. ---
\newblock Т.~22, {№}~5. ---
\newblock {С.}~663--670.

\bibitem{Bojanov75}
\selectlanguageifdefined{russian}
\BibEmph{Боянов~Б.~Д.} Наилучшие методы
  интерполирования для некоторых классов
  дифференцируемых функций~/ Б.~Д.~Боянов~//
  \BibEmph{Матем. заметки}. ---
\newblock 1975. ---
\newblock Т.~17, {№}~4. ---
\newblock {С.}~511--524.

\bibitem{Bojanov}
\selectlanguageifdefined{russian}
\BibEmph{Боянов~Б.~Д.} Оптимальное
  восстановление дифференцируемых
  функций~/ Б.~Д.~Боянов~// \BibEmph{Матем. сб}. ---
\newblock 1990. ---
\newblock Т. 181, {№}~3. ---
\newblock {С.}~334--353.

\bibitem{Zhensikbaev}
\selectlanguageifdefined{russian}
\BibEmph{Женсыкбаев~А.~А.} Приближение
  дифференцирумых периодических функций
  сплайнами по равномерному разбиению~/
  А.~А.~Женсыкбаев~// \BibEmph{Матем. заметки}. ---
\newblock 1973. ---
\newblock Т.~13, {№}~6. ---
\newblock {С.}~807--816.

\bibitem{ZhensikbaevBook}
\selectlanguageifdefined{russian}
\BibEmph{Женсыкбаев~А.~А.} Проблемы
  восстановления операторов~/
  А.~А.~Женсыкбаев. ---
\newblock Ижевск~: Институт компьютерных
  исследований, 2003. ---
\newblock 412~{с.}

\bibitem{hahn-banach}
\selectlanguageifdefined{russian}
\BibEmph{Колмогоров~А.~Н.} Элементы теории
  функций и функционального анализа~/
  А.~Н.~Колмогоров, С.~В.~Фомин. ---
\newblock М.~: Наука, 1976. ---
\newblock 543~{с.}

\bibitem{Kolmogorov85}
\selectlanguageifdefined{russian}
\BibEmph{Колмогоров~А.~Н.} О неравенствах между
  верхними гранями последовательных
  производных функции на бесконечном
  интервале~/ А.~Н.~Колмогоров~// Избр. тр.
  Матем, мех. ---
\newblock М.~: Наука, 1985. ---
\newblock {С.}~252--263.

\bibitem{Kornejchuk}
\selectlanguageifdefined{russian}
\BibEmph{Корнейчук~Н.~П.} Поперечники в $L_p$
  классов непрерывных и дифференцируемых
  функций и оптимальные методы кодирования
  и восстановления функций и их
  производных~/ Н.~П.~Корнейчук~// \BibEmph{Изв. АН
  СССР, Серия Матем}. ---
\newblock 1981. ---
\newblock Т.~45, {№}~2. ---
\newblock {С.}~266--290.

\bibitem{ExactConstants}
\selectlanguageifdefined{russian}
\BibEmph{Корнейчук~Н.~П.} Точные константы в
  теории приближения~/ Н.~П.~Корнейчук. ---
\newblock М.~: Наука, 1987. ---
\newblock 423~{с.}

\bibitem{MLD}
\selectlanguageifdefined{ukrainian}
\BibEmph{Моторный~В.~П.} Оптимальное
  восстановление функций и функционалов~/
  В.~П.~Моторный, А.~А.~Лигун, В.~Г.~Доронин. ---
\newblock Днiпропетровськ~: Вид-во ДДУ, 1994. ---
\newblock 224~{с.}

\bibitem{smolyk}
\selectlanguageifdefined{russian}
\BibEmph{Смоляк~С.~A.} Об оптимальном
  восстановление функций и функционалов на
  них~: {Дисс\ldots\ кандидата наук}~/
  С.~A.~Смоляк~; МГУ. ---
\newblock Mосква, 1965.

\bibitem{Tikhomirov87}
\selectlanguageifdefined{russian}
\BibEmph{Тихомиров~В.~М.} Итоги науки и техн~/
  В.~М.~Тихомиров~// \BibEmph{Сер. Соврем. пробл.
  мат. Фундам. направления}. ---
\newblock 1987. ---
\newblock Т.~14. ---
\newblock {С.}~103–--260.

\bibitem{Tihomirov}
\selectlanguageifdefined{russian}
\BibEmph{Тихомиров~В.~М.} Некоторые вопросы
  теории приближений~/ В.~М.~Тихомиров. ---
\newblock М.~: Изд-во МГУ, 1975. ---
\newblock 304~{с.}

\bibitem{borsuk}
\selectlanguageifdefined{english}
\BibEmph{Borsuk~K.} Drei S{\"a}tze {\"u}ber die n-dimensionale euklidische
  Sph{\"a}re~/ K.~Borsuk~// \BibEmph{Fund. Math.} ---
\newblock 1933. ---
\newblock Vol.~20. ---
\newblock P.~177--190.

\bibitem{Ligun}
\selectlanguageifdefined{russian}
\BibEmph{Ligun~A.~A.} Inequalities for upper bounds of functionals~/
  A.~A.~Ligun~// \BibEmph{Analysis Math.} ---
\newblock 1976. ---
\newblock Т.~2, {№}~1. ---
\newblock {С.}~11--40.

\bibitem{Micchelli76a}
\selectlanguageifdefined{english}
\BibEmph{Micchelli~C.~A.} On n-Widths and Optimal Recovery in $M^r$~/
  C.~A.~Micchelli, A.~Pinkus~// Approximation Theory, II~/ Ed.\ by\
  G.~G.~Lorentz, C.~K.~Chui, L.~L.~Schumaker. ---
\newblock N. Y., 1976. ---
\newblock P.~475--478.

\bibitem{Micchelli76b}
\selectlanguageifdefined{english}
\BibEmph{Micchelli~C.~A.} The optimal recovery of smooth functions~/
  C.~A.~Micchelli, T.~J.~Rivlin, S.~Winograd~// \BibEmph{Numer. Math.} ---
\newblock 1975. ---
\newblock Vol.~26, no.~2. ---
\newblock P.~191--200.

\bibitem{Micchelli77}
\selectlanguageifdefined{english}
\BibEmph{Micchelli~C.~A.} Optimal Estimation in Approximation Theory~/
  C.~A.~Micchelli, T.~J.~Rivlin. ---
\newblock NY~: Plenum Press, 1977. ---
\newblock 300~p.

\bibitem{Osipenko}
\selectlanguageifdefined{english}
\BibEmph{Osipenko~K.~Yu.} Optimal Recovery of Analytic Functions~/
  K.~Yu.~Osipenko. ---
\newblock Huntington, New York~: Nova Science Publishers, Inc., 2000. ---
\newblock 220~p.

\bibitem{Pinkus}
\selectlanguageifdefined{english}
\BibEmph{Pinkus~A.} N-widths and optimal recovery~/ A.~Pinkus~//
  \BibEmph{Proceedings of Symposia in Applied Mathematics}. ---
\newblock 1986. ---
\newblock Vol.~36. ---
\newblock P.~51--66.

\bibitem{Traub}
\selectlanguageifdefined{english}
\BibEmph{Traub~J.~F.} A general theory of optimal algorithms~/ J.~F.~Traub,
  H.~Wo\'{z}niakowski. ---
\newblock Academic Press, 1980. ---
\newblock 341~p.

\end{thebibliography}

\end{document}